\documentstyle[fullpage, twocolumn, twoside, psfig]{article}
\newtheorem{defi}{Definition}[section]

\newtheorem{prop}[defi]{Proposition}

\newtheorem{exa}[defi]{Example}
 
\newtheorem{rem}[defi]{Remark}
\newtheorem{question}[defi]{Question}

\def\N{\mathbf{N}}
\def\Z{\mathbf{Z}}
\def\Q{\mathbf{Q}}
\newcommand{\A}{\mathcal{A}}
\newcommand{\C}{{\mathcal{C}}}
\newcommand{\M}{{\mathcal{M}}}
\newcommand{\E}{\mbox{E}}
\newcommand{\code}{\mbox{\textsf{Code}}}
\newcommand{\codeA}{{a}}
\newcommand{\codeB}{{b}}
\newcommand{\inter}{\mbox{\textsf{Int}}}
\newcommand{\obs}{\mbox{\textsf{Obs}}}
\newcommand{\w}{w}
\newcommand{\uu}{{\mathbf{u}}}
\newcommand{\mufat}{{\mbox{\boldmath$\mu$}}}
\newcommand{\wfat}{{\mathbf{w}}}
\newcommand{\vv}{{\mathbf{v}}}
\begin{document} 
\title{In search of an evolutionary coding style}
\author{Torbj\"orn Lundh\thanks{%
   This work was done at SUNY Stony Brook under the 
   support (408003356-8) from the Swedish Natural Science Research Council, 
   NFR.}\\
 Department of Mathematics\\
 S-412 96 G\"oteborg, Sweden \\ torbjrn@math.chalmers.se}
\date{\relax}
\maketitle
\thispagestyle{empty} \def\IMSmarkvadjust{0 pt}
\def\IMSmarkhadjust{0 pt}
\def\IMSmarkhpadding{0 pt}
\def\IMSpubltext{Published in modified form:}
\def\SBIMSMark#1#2#3{
 \font\SBF=cmss10 at 10 true pt
 \font\SBI=cmssi10 at 10 true pt
 \setbox0=\hbox{\SBF \hbox to \IMSmarkhpadding{\relax}
                Stony Brook IMS Preprint \##1}
 \setbox2=\hbox to \wd0{\hfil \SBI #2}
 \setbox4=\hbox to \wd0{\hfil \SBI #3}
 \setbox6=\hbox to \wd0{\hss
             \vbox{\hsize=\wd0 \parskip=0pt \baselineskip=10 true pt
                   \copy0 \break%
                   \copy2 \break%
                   \copy4 \break}}
 \dimen0=\ht6   \advance\dimen0 by \vsize \advance\dimen0 by 8 true pt
                \advance\dimen0 by -\pagetotal
	        \advance\dimen0 by \IMSmarkvadjust
 \dimen2=\hsize \advance\dimen2 by .25 true in
	        \advance\dimen2 by \IMSmarkhadjust

%
%
  \openin2=publishd.tex
  \ifeof2\setbox0=\hbox to 0pt{}
  \else 
     \setbox0=\hbox to 3.1 true in{
                \vbox to \ht6{\hsize=3 true in \parskip=0pt  \noindent  
                {\SBI \IMSpubltext}\hfil\break
                \input publishd.tex 
                \vfill}}
  \fi
  \closein2
  \ht0=0pt \dp0=0pt
 \ht6=0pt \dp6=0pt
 \setbox8=\vbox to \dimen0{\vfill \hbox to \dimen2{\copy0 \hss \copy6}}
 \ht8=0pt \dp8=0pt \wd8=0pt
 \copy8
 \message{*** Stony Brook IMS Preprint #1, #2. #3 ***}
}
 
\def\IMSmarkhadjust{\hsize}  
\SBIMSMark{2000/3}{March 2000}{}


\begin{abstract}
%
%
In the near future, all the human genes will be identified. 
But understanding the functions coded in the genes is a much
harder problem. For example, by using block entropy, one has that the
DNA code is closer to a random code then written text, which in
turn is less ordered then an ordinary computer code; see \cite{schmitt}.

Instead of saying that the DNA is badly written, using our
programming standards, we might say that it is written in a 
different style --- an evolutionary style.

We will suggest a way to search for such a style in a quantified manner 
by using an artificial life program, and by giving a 
definition of general codes
and a definition of style for such codes.
\end{abstract}

\section{Background} \label{sec.background}
Let us as a background cite three different sources.

The first is  J. Madox's  comment on page 376 in \cite{madox}.
\begin{quote}
``In genetics for example, the task of understanding the functions of
all the 100,000 human genes will require a much greater effort than
that involved in their identification, and by a factor 10 or more.''
\end{quote}
Chris Adami from  Caltech made, in his
survey talk at Renaissance Technologies in Stony Brook 
10/27/98 about artificial life,  
a brief remark about the 
quality of the evolved program codes in his avida set up. 
He said something like this:
\begin{quote} 
``The codes that are evolved will eventually be almost totally unreadable. 
Things are never used only once, but two or more times. 
It is a kind of a `madman's' code.''
\end{quote}

In recent years, a lot of examples have been found where genes have been
``reused'' for different purposes during development. As an example,
take  the {\em runt} gene
in Drosophila, which is used in 
sex determination, segmentation and central nervous system creation.
In \cite{runt} the authors write:
\begin{quote}
``As mature organisms we are composed of an astonishing array of
diverse cell types---all derived from a single-celled zygote. When
faced with the task of generating such cellular diversity in a
reproducible fashion, how has the embryo chosen to respond? Recent
work in a number of developmental systems has suggested that the
embryo has employed two approaches. First, given finite resources, the
embryo has efficiently chosen to reutilize a limited set of proteins
in different temporal and spatial contexts to create cellular
diversity. Second, the embryo has also chosen to install molecular
redundancies to ensure the reproducibility of these patterns from
individual to individual.''
\end{quote}

How can we capture these comments about the style or quality of the
computer code and the DNA, in a quantified manner? Can we do that in
such a general manner that we will be able to use analogous quality
measure both for carbon-- and silicon based genetic codes?

\subsection{Plan of the paper} We  give a straightforward 
but general definition of a ``code'', a general definition of the
``style'' of such codes using a given set of measures. We will
also suggest such measures of
characteristic features, such as some type of 
``madness'', the robustness, and
the amount of reuse, of such codes.
Finally, we will make some ``baby--experiments'' by running the
avida program and  analyze samples of the evolved code,
generate two simulated programs in the same function class,
and eventually stylistically compare the real evolved code with these
simulated codes, all of this will be done in a C$++$ program. 
Some results will be graphically displayed at the
end.
\subsection{Future goals}
One could then more systematically  run the
avida program (or something similar) and 
analyze the evolved codes stylistically
with more realistic, and a higher number of  comparasion codes. By
changing parameters for the set up in avida, one might eventually
capture some common features.
By comparing the carbon based programming style and the silicon based
style, it might be possible to find some common parts that would
describe the natural programming style for evolution. That would in
turn help us read carbon based code. And it would  also give us
some hints how to create more robust computer programs, by looking at
how nature has solved such problems.

\subsection{Questions} Is there an existing way to 
study the style or quality of
the DNA? And is there any existing theory in the information sciences
that deals with this on the silicon side? We are aiming at a
complexity level higher than the usual information theory measures, such as
the notion of entropy (see the comment in Section \ref{sec.dna}).
\subsection{Acknowledgement}
I would very much like to thank professor G.\ Thomsen for 
listen to my ideas and for his constant encouragement. I would
also like to thank him for making me feel so much at home in his
{\em Xenopus} lab, and teaching me the basics in frog handling and
 cultural behaviour in the life sciences.

I would also like to thank Dr D.\ Slice for his interest and
all his help and suggestions.

B. Cohen has been a good source for me concerning the existing
computer code complexity measures. He has also been a 
fruitful discussion partner.

I am  grateful to S. Sutherland for reading an earlier version of this paper
and giving me  useful comments and pointing out errors and weaknesses.

Finally I would like to thank D. Brander for playing squash with
me and for trying to teach me some grammar. (All the remaining
language errors are naturally completely my own fault.)

\section{A general code}
We will give a general definition of a code as a string of generalized
``letters'' from a given alphabet, that when interpreted, will 
define a function. This interpretation is not unique, i.e.\ 
many different codes will produce the same function.
That fact will give us a way to study classes of codes.
That is, two codes are in the same class if their interpretation
gives equivalent functions.

\subsection{Codes}
We will give a definition of a code by using an underlying 
{\em alphabet}, and a generalized {\em interpretor}.
\begin{defi} \label{def.code}
Let us define a {\bf code} to be a finite string of ``letters'' taken from an 
``alphabet'', $\A$,
\[\code=\{\alpha_i\}_{i=1}^k, \mbox{ where } \alpha_i \in \A,\]
such that when the code is  interpreted, the code will represent a well 
defined
 function (or a process), $\code_j \rightarrow f_j$, with
a domain $D_{f_j}$ such that for all inputs, 
$x \in   D_{f_j}$, to the interpretation
of the code, will give $f_j(x)$ as the output\footnote{Note that
we here consider a function in a general sense, i.e.\ not necessarily
numerical.}.
\end{defi}

\subsection{Classes of codes}
From the above definition of codes, we see that different codes
can have the same function representation, see Example \ref{ex.fg}
below.

Let us therefore introduce the following classification.
Let the code $\code_f$  and the codes $\code_{f_i}$ have the functions
$f$ and $f_i$ as their representations.
\begin{defi} \label{def.c}
Let the code class with respect to the function $f$ be
the following (infinite) set of codes.
\[ 
\C_f = \{ \code_{f_i} : f_i(x)=f(x), \mbox{ for all } x \in D_f \}.
\]
\end{defi}

\begin{rem}
If the interpretor is not able to produce a well defined 
function from the code $\code$, we say that
$\code$ is in the error class $\C_\varepsilon$. That (huge) class
can be viewed as the complement to all interpretable codes. 
\end{rem}

\begin{question}
With the suggested norm of $\C_f$ below, is  $\C_f$  compact?
Do there exist extremal codes?
\end{question}

\begin{exa} \label{ex.fg}
Let 
\[f(x) = \frac{1}{x-2} ,\]
\[g(x)= \left\{ \begin{array}{ll}
                \frac{x+2}{x^2-4} & \mbox{ if $x\neq -2$}\\
                -\frac{1}{4} & \mbox{ if $x=-2$}
                \end{array}
        \right. \] and
\[h(x)= \left\{ \begin{array}{rlr}
                -\frac{1}{4} & \mbox{ if $x=$}&-2\\
                -\frac{1}{2} & \mbox{ if $x=$}&0\\
                -1 & \mbox{ if $x=$}&1\\
                1 & \mbox{ if $x=$}&3\\
                \frac{1}{2} & \mbox{ if $x=$}&4
                \end{array}
        \right. \]
Note that $\code_g$ and $\code_f$ are both in 
$\C_h$; and  $\C_f=\C_g \subset \C_h
$.
\end{exa}

\begin{rem}
One common way to view functions is as black boxes. Here we are interested
of the internal structure of such black boxes performing 
equivalent tasks.
\end{rem}

\section{The style of a code} \label{sec.style}
We will now give a  description of a method to characterize 
different coding styles.
Since different codes do not differ in function in the 
same class, we say that they differ in {\em style}. How can we
characterize such style?

Given a set of  measures on the codes, we 
propose a way to characterize style of a subset of codes
in the code class as an extremal unit weight vector that will act as
a stylistic  ``fingerprint''.
We will also give an algorithm for ``stylistic translations'', and
an index which reveals how well a given subset represents a 
common style.

The methods we use apply very elementary mathematics, and maybe even more
basic statistical methods. 
That will hopefully make it accessible to a wide scientific audience
who are interested in ``style''.

\subsection{A measure on $\C_f$}
Let us now study measures on  $\C_f$. 
Let 
\[\mu_i: \C_f \rightarrow [0,1].\]
Let us consider the following {\em profile measure}
\[\mufat=(\mu_1,\mu_2, \ldots, \mu_n)\]
 of  codes in  $ \C_f$.   
If a given measure has a range outside $[0,1]$, let us 
use the transformation
\[x \rightarrow \frac{x}{1+x}\]
to make it fit into $[0,1]$.

We define the following scalar measure.
\[\nu_{\wfat}(\code_g)= \wfat \cdot \mufat(\code_g),\]
where $\wfat$ is a normalized  weight vector such that
$||\wfat||=1$,
for some norm $||\cdot ||$.
For example, let
\[||\wfat ||=||\wfat||_p=\biggl( \sum_{i=1}^n |\w_i|^p\biggr)^{\frac{1}{p}},\]
where $p\geq 1$. 
As a default norm, let us use $||\cdot||=||\cdot||_2$.
 
\subsection{Extremal weights}
Let us use the above measure to try to capture a characterization
of ``style'' of codes.

Suppose that we have two sample sets, $A$ and $B$,  of codes  in $\C_f$.
We will try to 
find a style characterization of the codes in $A$ relative to  $B$.

We can think of $A$ as the set of codes we are stylistically
interested in and $B$ as a complementary environment.

Let us define a vector $\uu=(u_1,u_2, \ldots, u_n)$ in the 
following way.
\begin{equation} 
\label{eq.u} 
\uu=\sum_{\codeA_i\in A} \sum_{\codeB_j \in B} 
(\mufat(\codeA_i)-\mufat(\codeB_j)).
\end{equation}

Let us now normalize $u$ to a unit vector.
\begin{equation} \label{eq.w}
\wfat^+=\wfat^+(A)= \frac{\uu}{||\uu||}.\end{equation}

Let us now study the random variable
\begin{equation} \label{eq.x} 
X=\nu_{\wfat}(\codeA_i)-\nu_{\wfat}(\codeB_j),
\end{equation}
where $\codeA_i$ is a randomly chosen code in $A$, 
with uniform 
probability\footnote{The probability to choose  $\codeA_i$ is $1/\#(A)$.}, 
and
 $\codeB_j$ is  randomly chosen  in $B$.
Note that $X$ is dependent on the chosen $\wfat$.

\begin{prop} \label{prop.w}
Let $X$ be the random variable defined in (\ref{eq.x}) and let
$\wfat^+$ be the unit vector from (\ref{eq.w}). Then we have
that picking $\wfat=\wfat^+$ will maximize the expected value of
$X$, $\E(X)$.
\end{prop}

{\bf Proof:}
Let $\#(\cdot)$ denote the cardinality of a set and let
$M=\#(A) \#(B)$.
The  expectation for a general weight vector $\wfat$ will be easily computed.
\[E(X)=\frac{1}{M}\sum_{\codeA_i\in A} \sum_{\codeB_j \in B} 
(\nu_{\wfat}(\codeA_i)-\nu_{\wfat}(\codeB_j))=\]
\[=\frac{1}{M}\sum_{\codeA_i\in A} \sum_{\codeB_j \in B} 
(\wfat\cdot\mufat(\codeA_i)-\wfat\cdot\mufat(\codeB_j))=\frac{\wfat\cdot \uu}{M}.\]
Since $||\wfat||=1$ we have that 
\[|E(X)|\leq\frac{||\uu||}{M}.\]
Let us now look at the special case $\wfat=\wfat^+$.
We have then that
\[E(X)=\frac{\wfat^+\cdot \uu}{M}=\frac{1}{M}\frac{\uu}{||\uu||}\cdot \uu=\]
\[=
\frac{||\uu||^2}{M||\uu||}=\frac{||\uu||}{M}.\]
Thus we see that taking $\wfat$ to be $\wfat^+$ will maximize $E(X)$. $\Box$


\subsection{A fingerprint on the hyper sphere}
We can now view the  unit vector $\wfat^+(A)=\wfat^+$ 
on the unit hyper\footnote{if $n>3$} sphere as
a ``fingerprint'' of the style of codes in $A$ relative to $B$ with
respect to the list of given measures in $\mufat$. 

\begin{rem}
If there is a superset of codes $B_1\supset B$ which is 
enough separated, then the fingerprint $\wfat^+$ of $A$
relative to $B_1$ can be expected to be a 
more refined characterization of the ``style''
in $A$ than the $\wfat^+$ of $A$ relative to $B$.
\end{rem}

\subsection{A universal character}
One would hope that similar classes, $\C_f$, would
give similar fingerprints $\wfat^+$ for different function codes created by 
the same code writing agent.

\subsection{Stability with respect to the vector $\mufat$}
There is a stability feature built into  $\wfat^+$ in the sense that
if you would like to find a stylistic quality in a group of codes,
you  try to find ``relevant'' measures in the vector $\mufat$. 
What happens if, in addition to your relevant measures, you also
        take a sequence of irrelevant measures (where, for example, the
        codes look more or less randomly distributed, or even just the
        same)? 
If the environment, i.e.\ $B$,
is rich enough, then your profile will just be about zero at the tail,
where the non-relevant measures are. That means you don't have to 
be restrictive when you pick your measures --- if some happen to 
be worthless, that will be taken care of by itself.

\subsection{Is there a common style in $A$?} \label{sec.common}
Given two subsets, $A$ and $B$ in $\C_f$, of codes, we have now a method for
finding a common style in $A$, in relation to $B$, as a unit vector $\wfat^+(A)$.
(As a special case, we can let $B$ be the complement of $A$ in $\C_f$.
Then we can talk about {\em the} style of $A$.)

This process  can  be executed even if $A$ and $B$ are 
just  randomly chosen subsets in 
$C_f$ where we can not expect to find any stylistic common features
in $A$ in comparison to all the codes in $B$. How can we 
find out if $A$ really has a common style, in comparasion to 
$B$, that can be captured
by the chosen measure profile $\mufat$?

Let us go back to Equation (\ref{eq.x}) and the random  
variable $X$. Let us also define a similar 
s.v.\ $Y$ such that
\[Y=\nu_{\wfat}(c_i)-\nu_{\wfat}(c_j),\]
where $c_i$ and $c_j$ are randomly chosen, with equal 
probability, in $A\cup B$. Note that the underlying $\wfat$ is 
$\wfat^+(A)$. From Proposition \ref{prop.w} above, the expected value of 
$X$ will then be maximal. Let us denote that value by $m$, i.e.\ let 
$ E(X)=m$. From the 
proof of Proposition \ref{prop.w} we see that \[m=\frac{||\uu||}{\#A \#B}.\]

Hence, loosely speaking, $m$ is large if $A$ has a characteristic
style that is captured by the measure profile $\mufat$. 
How large can $m$ get? Or in other words: how large can $||\uu||$ get?
Since all the measures $\mu_k(\cdot)$ are bounded above by $1$ and 
below by $0$ we have that a component $u_i$ of the $\uu$ vector
is also bounded.
\[u_k=\sum_{\codeA_i\in A} \sum_{\codeB_j \in B} (\mu_k(\codeA_i)
-\mu_k(\codeB_j))\leq \]
\[ \leq \sum_{\codeA_i\in A} \sum_{\codeB_j \in B} (1-0)=\#A\#B=M.\]
If we are using the standard norm $||\cdot||=||\cdot||_2$,
we get that
$||\uu||\leq n M^2$ and hence $||\uu||\leq \sqrt{n} M$.
  
Let us normalize $m$ to get 
the following index \[\theta=\frac{m}{\sqrt{n}}.\]

We see that $\theta \in [0,1]$   and it will be closer to $1$ when the 
fingerprint is good.

On the other hand, due to the  
Remark \ref{rem.k} below, we might also need an invariant index.

Let us use the 
variances of the random variables $X, Y$ in the following way. 
Let 
\[\sigma_A^2 = E((X-m)^2), \mbox{ and let } \sigma_{AB}^2= E(Y^2).\]
Note that we immediately have that $E(Y)=0$.
Let us suggest an index $\eta$ of how well the measure profile
 $\mufat$ can capture a common style, if it exists at all, 
in $A$ in comparison to  the codes in $B$.
Let 
\begin{equation} \label{eq.eta}
 \eta=\eta(A, B, \mufat)=\frac{\sigma_{AB}^2}{\sigma_{A}^2}=
\frac{E(Y^2)}{E(X^2)-m^2}.
\end{equation}

\begin{rem} \label{rem.k}
Note that if all the measures in $\mufat$ are multiplied by a factor $k<1$,
i.e. \ $\mu_i \rightarrow k\mu_i$, then $m \rightarrow km$,
$\sigma_A^2\rightarrow k^2\sigma_A^2$, and  
$\sigma_{AB}^2\rightarrow k^2\sigma_{AB}^2$. Hence
$\eta$ would not change, but $\theta \rightarrow k \theta$.
\end{rem}

\subsubsection{Principal component analysis}
A very illustrative, and very popular, technique when looking
for connections in a multidimensional environment is
the {\em principal component analysis}.
It is  a two dimensional diagram with coordinate
axis the two eigenvectors with the largest eigenvalues of the
covariance matrix $\M$. Let $C=A \cup B$ and let $N=\#C$. We
denote the codes in $C$  by $c_j$, then
\[\M_{i,j}=
\mbox{Cov}\biggl((\mu_i(c_1),\mu_i(c_2),\ldots, \mu_i(c_N)),\]
\[(\mu_j(c_1),\mu_j(c_2),\ldots, \mu_j(c_N))\biggr).\]

For more details, see any book in multivariate
analysis, for example \cite{murtagh}.

If in such a diagram, the codes $A$ we are interested in are
clearly separated from the $B$ codes, then we could say that
there is a common style in $A$, and if that is the case and
if the largest eigenvalue is considerably larger then the
second one, then the fingerprint $\wfat$ and the 
first eigenvector should be close to each other. That is indeed
the case in Figure \ref{fig.principal}.

\subsubsection{Cluster analysis}
Another tool from the multivariate toolbox could be used to 
study the question about a common style in $A$. This method
is based on an iteration fusion of close points until the
desired number of subsets are obtained. 
To measure the closeness, we might pick our scalar measure
$\nu_{\wfat^+}$. 

In order to check if there is a common style in $A$, we can ask how well
the points in $A$ are clustered.

\subsection{Style translations by iterations}
Using the above construction, let 
us here indicate an algorithm how to ``stylistically 
translate'' a code $\codeA$
in $A\subset \C_f$, into the style in $B \subset \C_f$.

\begin{enumerate}
\item Compute 
\[\vv=\sum_{\codeB_i \in B}  
(\mufat(\codeB_i)-\mufat(\codeA)).\]
\item Find the dominating component of the vector $\vv$, i.e.\ let
\[v_m=\max_{1\leq i\leq n} |v_i|.\]
\item Study the measure $\mu_m$ and stylistically rewrite $\codeA$ such
that $\mu_m(\codeA)$ would increase approximately $v_m$ units (decrease if
$v_m<0$).
\item Iterate the process until sufficient accuracy is attained.
The accuracy is measured by $||\vv||$.
\end{enumerate}

Suppose the final accuracy is $\delta$ and
denote the rewritten $\codeA$ by $\codeA'$. Let $\wfat=\wfat^+(B)$
and let $Z$ be the random variable 
\[Z=\nu_{\wfat}(\codeB_i)-\nu_{\wfat}(\codeA'),\]
for a randomly chosen $\codeB_i$ in $B$.
Now the expected value of $Z$ will  be
\[E(Z)=\frac{\sum_{i\in B}
\wfat \cdot (\mufat(\codeB_i)-\mufat(\codeA'))}{\#(B)},\]
where $\#(\cdot)$ stands for the number of elements in the set.
That is
\[ E(Z)=\frac{\wfat \cdot \vv}{\#(B)}.\]
Now if we have the norm $||\cdot||$ as the default
$||\cdot||_2$ we can use Cauchy--Schwarz' inequality to 
get that 
\[|E(Z)|\leq \frac{\delta}{\#(B)}.\]
In other words, the fingerprint $\wfat^+$ of $B$ would hardly
``feel'' the difference between $\codeA'$ and the codes in
$B$ if $\delta$ is small.

\section{Different levels of Code}
So far we have just been studying a code on a singular level. In
this section we will describe a way to separate the codes into
different levels. One could then ask questions about the styles on
different levels. Are the styles similar even on different levels, etc?
 
We will view a code as a  composition of 
lower level codes. Let us use the definition of
emergence given in \cite{baas} on p. 518.
\[ P \mbox{ is an emergent property of } S^2 \]
\[ \Longleftrightarrow \]
\[P \in \obs^2(S^2), \mbox{ but } P \not \in \obs^2(S^1_{i_1}) \;\;\;
\forall i_1.
\]
Let us explain this more in detail, and similarly 
display a concrete example
of a computer code, where the inputs are integers.

In this case, 
let $S^0$ be input values, e.g.\ $S^0_i\in \N$, 
and let $\{\inter^i\}_{i=0}^n$ be a given
sequence of sets of interactions, 
e.g.\ $\inter^0=\{+,-\}$, $\inter^1=\{*,\div\}\cup\inter^0$,
 $\inter^2=\{=,<,>\}\cup\inter^1$, etc.
Let  $\{\obs^i\}_{i=0}^n$ be the related observational function,
e.g.\ $\obs^0(x)=$ value of $x$ as an integer,  
$\obs^2(x>y)=$ true or false.

From the given sequences  $\{\inter^i\}_{i=0}^n$ and 
$\{\obs^i\}_{i=0}^n$ we get the higher order structures as 
a ``reaction'', $R$, in the follwing manner.
\[S^1=R(S^0, \obs^0, \inter^0), \]
\[S^2=R(S^1, \obs^1, \inter^1), \]
see \cite{baas} for details.

We have then that in our example, $-1$ is an emergent property of $S^1=\Z$, 
$1/2$ is an emergent property of $S^2=\Q$, and {\em ``true''} is an emergent
property of $S^3$, etc.

Given such sequence $\{\inter^i\}_{i=0}^n$ we can define
a code of degree $k$ as a consecutive string of the total code which
has the property of $\obs^k$, i.e. if the code can produce an
output that is an
emergent property of $S^k$. 

We can then view the final code  of degree $n$, $\code_n$, as a composition
of sub-codes of degree $n-1$, $\code_{n-1}^i,$ etc.

Note that Baas indicates this application in mentioning 
the word {\em hyperalgorithms} on p. 526 in 
\cite{baas}.

How can one think of a good implementation of these functions?
The choice of the interactions $\inter$ will give us a chance
to find fine structures. What happens to our example if we 
start by $\inter^0=\{+\}$ and then $\inter^1=\{-\}\cup\inter^0$ etc?
Is there a ``natural'' choice of that interaction sequence
for a given case?


\section{Some measures}
Let us list a couple of measures that would be useful to capture
some of the features mentioned in the 
Section \ref{sec.background} above, and which also utilizes
different levels of codes discussed above.

\subsection{Spaghetti}  \label{sec.spatti}
Let us propose a kind of ``code-madness-measure'' using the 
above hierarchies of codes.

Let 
$m_k$ be the maximum numbers of $\code_{k-1}$ codes in 
a $\code_k$. More precisely,
Let 
\[\eta^i_k=\max \{ j: \code_{k-1}^{i,j} \subset \code_k^i \},\]
and let  
\[m_k= \max_i \eta_k^i \mbox{ and } s_k= \sum_i \eta_k^i.\]
Now, we define the ``spaghetti length'' at level $k$ to be
\[S_k(\code_n)=\frac{m_k}{s_k} \mbox{ and } \]
\[S(\code_n) = 
\max_k S_k(\code_n).\]

\subsection{Reuse}
Let $0<k\leq n$
\[R^i_k(\code_n)=\]
\[\frac{\mbox{max   $\# \code_{k-1} \subset \code_k$   
used $i$ times or more}}{s_k}.\]

Let us  use the conventions
$R^2_k(\code)=R_k(\code)$ and $R^i_2(\code)=R^i(\code)$.

\subsection{Redundancy}
A very important feature in evolutionary driven codes are
their robustness.
There are good measures of robustness given in the literature.
One way to measure it is to to check the probability that
the code ``survives'' a one point mutation.

The redundancy in a code of degree $k$,
could be measured in the following way.

Let $m$ be the maximal number of subunit codes of type $\code_{k-1}$
in the code of type $\code_k$, that can be taken away without affecting
the output of the $\code_k$ code, and let $n$ be the
total number of subunits. Let us then define the 
redundancy in  code $C$ of type $\code_k$ by
\[Red(C)=\frac{m}{n}.\]

\subsection{Brittleness}
Another function that might be better to use in some ways would be
a ``brittleness'' function that we define in the following way.

Let $m$ and $n$ be as above and let $d(k)$ be the number of 
$\code_{k-1}$ codes in a $\code_k$ code, that when removed totally, or
partially, destroys the  code $\code_k$. The we define the brittleness of 
$\code_k$ as.
\[\mbox{\sl Britt}(\code_k)=\frac{d(k)}{n-m}.\]
That can be viewed in the following way. Let us think of $\code_k$ as a chain
and the subunits as links in that chain. Then the length of the chain will
be $n-m$ and $d(k)$ will be the number of links  that are can not be removed
without breaking the chain.

Note also that in the special case when we either have one or two links
for each step in the chain, that $d(k)=n-2m$.

Note that  the probability to survive a deletion of a randomly chosen subcode
is $1-d(k)/n$.
That is, if there would just be one sublevel to the code,  $1-d(k)/n$
would be the usual robustness measure mentioned above.

\section{Five applications}
Since we are mainly interested in searching for a 
evolutionary universal coding style, we are of course
interested in the two special cases when the code is
a computer code, and when it is a sequence of the DNA. 
We have tried to make the above definitions of codes and
style general enough to be able to deal with those cases.

As a by-product, we noted that this theoretical framework
could also be applied to other areas where `style' is essential.
The applications that comes closes to our mind was art, music and
literature. As a third application we will comment how the above
stylistic fingerprint could be used for a common framework in
``stylometry'' in the study of authorship attribution. Our 
hope is that if a theory is not only applicable to the special
cases it was aimed for under its construction, but also to
a different case, it might be a sound approach in that theory.

A fourth possible application would be a stylistic investigation of
the internal architecture of black boxes in the theory of neural nets,
and its carbon based version---the brain. 

The fifth application is
about artificial life. We will make a small experiment in such an environment.


\subsection{Computer code}

By looking at the indicated toy-example above where we had natural numbers
as inputs and $+,-$ as the first interactions, it is not to hard to
imagine that you would get the intuitive ``usual meaning'' increase in
complexity in substructures in the code from 
lower order arithmetic operations, more complex functions, subroutines,
program parts, the complete program. 

Note that what usually is seen as a good programming style, that is
separate codes into small, more or less, independent units, and not
too many jumps back and forth, will give you a low $S(\code)$; see
Section \ref{sec.spatti}.

Let us now mention three existing families of computer code measures.
See also \cite{jones} and \cite{ency} for a  detailed 
description of the two first examples and many others.

\subsubsection{Halstead's Complexity Measures}
In 1977 M.\ Halstead introduced a tool to measure the
complexity of a computer program. It is  
perhaps the most well known
measure of that kind.

The measure, or more precisely the family of five measures, 
is based directly on the  code in the following
way. Let $n_1$ be the number of distinct operators, 
$n_2$ the number of distinct operands, $N_1$ the total number
of operators, and $N_2$ the total number of operands.
From these numbers, the following measures are 
constructed.

\vspace{3mm}
\begin{tabular}{ll}
Program vocabulary & $n=n_1+n_2$.\\
Program length & $N=N_1+N_2$.\\
Difficulty &$D=\frac{n_1 N_2}{2 n_2}$.\\
Volume & $V=N \log_2(n)$.\\
Effort &$E=DV$.\\
\end{tabular}
\vspace{3mm}

The Halstead's measure seems to be a good candidate for measures in
the stylistic search since it is only based on textual information and
due to the fact that it has been used in many contexts and over such
a long time and hence its properties are quite well known.

\subsubsection{McCabe's Cyclomatic Complexity}  
In 1976, T. McCabe introduced a measure
of the number 
of linearly independent
paths through a computer program as a measure of the complexity of the
code.

The cyclomatic complexity, $CC$, is computed in the following way.
Let us study a schematic graph of the program and count the
number of edges $E$, the number of nodes $N$, and the number of 
connected components $c$. Then
\[CC=E-N+c.\]
If a program has a $CC$ higher than $50$ it is said to be unstable, 
since 
 it is then ``very likely'' break down if it is 
altered. 

The Cyclomatic Complexity is more of a measure of the inner logical
complexity then the textual Halstead's measure.

\subsubsection{GRASP}
Let us describe an ongoing project at Auburn, Alabama, which address
a numerical local measure of computer codes, and its display. 
See \cite{grasp} for further information and down--loading of
the program.
\begin{quote}
``The overall goal of the GRASP project is to improve the
comprehensibility of software. Thus, it is important to be able to
identify complex areas of source code. The Complexity Profile Graph
(CPG), a new graphical representation based on a composite of
statement level complexity metrics provides the user with the
capability to quickly recognize complex areas of source code. The CPG
is significant in that it shows the complexity of a program unit as a
profile of statement level complexity metrics rather than as a single,
global metric.'' 
\end{quote}

First the program code is parsed into non overlapping segments; then
a series of measures, briefly described below, is applied to each
segment, giving a local complexity measure of the program code.

The content complexity of a segment $S$ in the code is defined as
\[\eta(S)= log(\sum_{T\in S} \mbox{Weight(T)}),\]
where $T$ are tokens in the segment $S$. For example, 
in \cite{grasp}, the weights for Ada 95 are 
given in the table below.
\vspace{.3cm}

\hspace{-.7cm} 
{\small
\begin{tabular}{|l|l|c|} \hline
 Token Description & Symbol & Weight \\ \hline
  Logical operators &   and, or, not, ...&     1.5 \\ \hline
          Comparison op.&  $<, >, =, <=$, ...&       1.5\\ \hline
          Left parenthesis&   (&   1.3\\ \hline
          Identifiers&   var1, proc1, ...&        1.0\\ \hline
          Others&  $+, -, \cdot, /, ), \ldots$& 1.0\\ \hline
\end{tabular}
}

\begin{quote} 
``The context complexity provides a baseline level of complexity for
segments of simple statements nested within a compound statement,
which itself may be nested several levels deep.
The complexity of a compound statement is based on three aspects:
inherent complexity, reachability, and breadth.''  
\end{quote}

These three
complexities are 
added, with weights, to obtain the context complexity.

Combining the content complexity and the context complexity, by a 
weighted sum,
gives the profile metric $\mufat(S)$ for a segment $S$.

The local metric $\mufat$ can then be graphically presented as
a histogram to give indications where the code is more 
complex. Thus that would help a 
programmer to point out code segments that would need some extra
thought.

\subsubsection{Finding a coding style}
Suppose that we would like to find a way to assign a specific style to
a programmer. How could we do that? 

Suppose he writes in C++. Then
one might take a sample of his codes, look at his functions and 
subroutines. Go to other programmers and find as many functions as 
possible from that environment that one would like to be able
to identify ``our'' programmer in the future. 

Sort the codes into function classes $\C_f$ where our programmer has
written $\codeA_{f} \in \C_f$ and other programmers  $\codeB_{f,j} 
\in B \subset \C_f$.

Let us now take a wide variety of measures of computer codes, such as those
mentioned above, e.g.\ Halstead's measures etc. Let us call the 
vector of those measures $\mufat$, where each $\mu_i$ is a  
function of a C++ code to the unit interval $[0,1]$.

Let use  equation (\ref{eq.u}).
\[\uu=\sum_{f} \sum_{\codeB_{f,j} \in B} 
 (\mufat(\codeA_{f})-\mufat(\codeB_{f,j})).\]

Now, 
\[\wfat^+=\frac{\uu}{||\uu||}\]
will be the stylistic fingerprint of this programmer with respect to 
the vector of measures $\mufat$ and in the chosen environment. Note that
$\wfat^+$ will be heavily dependent on the chosen reference codes,
so it is important that that environment is chosen carefully.

\subsection{DNA} \label{sec.dna}
It is hard to get a handle on the DNA code with usual computer
code tools. As an example in \cite{schmitt} the authors 
study {\em block entropies} 
for DNA to find some indication of structure.
\begin{defi}
Let the length of the alphabet be $\lambda$. Then the 
(normalized) $n$--block entropy of a sequence is defined as
\[H_n=- \sum_{i=1}^{\lambda^n} p_i^{(n)} \log_\lambda p_i^{(n)} ,\]
where $p_i^{(n)}$ is the probability of the $i$th combination of $n$ ``letters''. 
\end{defi}
 In
the summary of \cite{schmitt} one can read the following.
\begin{quote}
``Surprisingly, DNA sequences behave closer to completely random sequences
than to written text. The very strict syntax of computer languages on the
other hand is reflected by a very low average information content of its
sub-strings.''
\end{quote}

To us that means that the style of DNA sequences is not just close to 
chaos, but rather written in such a different style compared to
written text, and even more different from computer codes.

The inputs, $S^0$ here are the amino acids, and one could take perhaps
$\inter^0$ to be ``putting next to'', or sequencing.

The letters on the DNA level should be the four bases, the
words are then three letters word coding for a amino acid.
Add to the dictionary special start and stop words for the
genes.

One could also study the situation on the amino acid level, i.e.\ one
level up from the bare DNA code. That would give us 20 letters in an
alphabet.\footnote{The alphabet is not unique, e.g. different spelling
in example some bacterias and humans for some amino acids.} The 
words would here be the genes.

\subsubsection{Hopes}
One would hope that some hypothetical insight of a universal 
evolutionary driven style would give some hints how to 
better ``read'' the DNA code. 

We would therefore try to find measures that are general enough in the
below described computer experiment, so that some hypothesis about
the DNA code could eventually be made. 

There are indications that evolutionary driven code, may not be
``optimal'' in a basic sense, but might include peculiar turns and
twists.
See for example p.\ 180 in \cite{madox} where
J.  Madox discuss {\em RNA Editing}
\begin{quote}
``The puzzle is to know why these changes, which are presumably
advantageous to the organism, have not been incorporated in the 
gene themselves, thus avoiding the need for editing by way of 
afterthought---not to mention the need for a separate biochemical
mechanism for carrying it out.''
\end{quote}
He also addresses, on p. 203,  the ``junk code'' in eukaryotic cells.
\begin{quote}
``At the very least, this complication is an extra metabolic cost for 
eukaryotic cells. It is also potentially a source of error. What 
countervailing selective advantage can there possibly be in this
arrangement?''
\end{quote}

\subsection{Literature}
Let us now turn our attention to something different. Let
us look at some examples in literature where ``style'' has been
in the focus.

On p. 74 in \cite{ravin} the author discuss computational differences
between grammatical errors and stylistic weaknesses.

In \cite{rudman}, seven important problems with the existing 
authorship attribution studies are listed and discussed and
some solutions proposed. The proposed solution to problem number
three is to 
\begin{quote}
``study style in its totality. Approximately 1,000 style markers
have already been isolated. We must strive to identify all of the markers
that make up ``style'' --- to map style the way biologists are mapping
the genes.''
\end{quote}
Furthermore, the suggested solution to problem five is to 
\begin{quotation}
``Develop a complete and necessarily multi--faceted theoreti\-cal
frame\-work on which to hang all non-traditional authorship attribution
studies.  

Publish the theories, discuss the theories, and put the theories to 
experimental tests.''
\end{quotation}

As an example of a suggested metric from the literature studies, one
can take the Yule's coefficient advocated among many others in \cite{delcourt}.
Let $\{f_{ij}\}$ represent the observed frequencies in a 
``two way contingency table'' and let
the Yule's coefficient (see \cite{yule}) be defined as
\[Y=\frac{\sqrt{c-1}}{\sqrt{c+1}}, \mbox{ where } 
c=\frac{f_{11} f_{22}}{f_{12} f_{21}}.\]

\subsubsection{Stylometry: Finding the author}
Let us now treat a hypothetical case, using the methods from Section 
\ref{sec.style}, on the authorship attribution problem.

What do we mean by a code in this case? Let us look at Definition 
\ref{def.code}. We will interpret the ``letters'' as word taken
from the ``alphabet'' which will in this case be a complete dictionary of
the language in question. What about the functions? 

Scarry gives in \cite{scarry} a description how we can
think about the reading process; see the following 
quote from the first chapter:
\begin{quote}
``When we say
      'Emily Bront\"e describes Catherine's face,' we might also
      say 'Bront\"e gives us a set of instructions for how to imagine
      or construct Catherine's face.' This reformulation is
      accurate if cumbersome, in that it shifts the site of mimesis
      from the object to the mental act.'' 
\end{quote}

So, in this case we might think about the functions as descriptions,
or simply constant functions, 
where the interpreter is
the reader.\footnote{Here it is easy to see the crucial importance of
the interpreter. I would for example be an extremely poor
interpretor of a French text, even if I would get some 
picture of a story at the end.} Examples of interpreted functions
would be text where ``boy meets girl'', or ``prince meets ghost''.
An even more refined version would be ``boy meets girl 
described in
an English sonnet''.

Suppose now that we want to test the hypothesis that author X
has written a given sonnet $s$. Then we might gather all know 
sonnets written by X into the set $A$ as a subset of 
reference sonnets from that time period  in $S$.
Now $\codeA_i \in A$ means sonnet number $i$ in $A$.

One wants to have as large environment $S$ as possible, but at
the same time also as narrow as possible, i.e.\ from the same time period, etc.
This has then to be chosen carefully and with great knowledge about the
literature period and its authors. Let now $B$ be the
set of comparison sonnets.

On the other hand when it comes to choose metrics, then we pick as many as
it is numerically computable in reasonable time, which depends both on
our time and our computer. For example, 
there has been some new interesting development
using neural networks in authorship attributions, see for example
\cite{tweedie}. An even more exciting method was 
used in \cite{forsyth} where they used genetic algorithms 
in order to find the best features to measure. 
We can use such nets and genetically derived measures in  
our set of measures too!

Now, apply equation (\ref{eq.u}) to get the vector $\uu$ which then
is normalized to $\wfat^+$. That will be the stylistic finger
print of author X, given the above constructed framework.

How good and reliable is this fingerprint? Suppose that $A$ is
large and $S$ is rich, not only large but more or less complete with
respect to the author representation. Then $\eta$ from equation (\ref{eq.eta})
would be a good measure how reliable $\wfat^+$ is. We want 
$\eta$ to be large of course, but what is large enough?

To answer that we need to study the actual distributions of the
sonnets evaluated by
$\nu_{\wfat}$, make approximations and perform hypothesis tests.

The nice thing about this method is that one does not have
to argue which measure should be applied. Just apply them all!

Talking about measures, there has been a very animated debate about
the values of using computers in literary studies. R. Quiones said in
\cite{atlantic} ``Why don't they simply read the plays?''. A literature
expert reading a text
is undoubtly very hard to beat when it comes to  authorship
attributions. But let us look at this situation having our suggested
definition of measures of codes in mind. Isn't the expert using a
long array of measures, weighted together in an intricate and sometimes
subconscious way? One measure could be: ``X would never use that word in 
two consecutive sentences''.  The measure would be the characteristic
function of that event in the text, and the weight would be heavily negative.
The more skill and familiarity the expert has about such tasks, the
lager the array $\mufat$, and the more subtle the weighing process would be.
It may seem like a trivialization to think in those terms, but since
the task is very hard and complex, one would not expect that a computer
would be programmed in the near future that would be generally better
than a literature expert. Compare with the relatively straightforward 
problem to play chess.

\subsection{Neural nets}
P. Adams in \cite{adams} 
pointed out
the possibility of strong evolutionary forces acting inside the
brain in the process of learning. ``Good'' synapses will be rewarded
by being strengthened, and ``bad'' will be punished by being weakened.
What is good and bad are much more intricate and implicit qualities compared
to the genetic evolution, but maybe evolutionary forces are in command
in the learning process on a time scale of seconds instead of millenias.
And if there is a universal coding style---that should then be found in the
style of the programming of the neural nets too.

\subsection{Artificial life}
Avida, see \cite{avida} and \cite{adami}, 
is a program in the Tierra class, see \cite{tierra}. 
Unlike Tierra it
gives a natural 2 dimensional picture of an evolutionary 
process. Like Tierra, Avida is {\em not} a simulation
of real carbon based DNA, but a real evolution in
a silicon based (simple) world. We will now do an experiment in the Avida world.

\section{A small experiment}
Let us try some of the above concepts in
a Avida experiment. We will deliberately 
try to make it as simple as possible in
order to get some output in a straightforward
way. But we will keep all the doors open to
make variants and generalizations in the future.

We use Avida version 1.0.1 that is available on the web;
see \cite{avida}.

When running the program, you have the possibility to extract
individuals and saving their data and their genetic code. To 
read more about this, see the documentation on \cite{avida} or
even better---the book \cite{adami}. The code and the data is
saved to files such as \texttt{153-aagxs}. 

The data in 153-aagxs tells you, for example, that in addition to
        being able to replicate itself, the code also performs some other
        tasks. 
In this specific case, it takes input from a stack, performs a logical XOR
        twice, performs three NOTs, etc. 
For this, it is
rewarded with a bigger time slice, and hence will reproduce (and
survive) better.

\subsection{The function class}
Let us as the function $f$ defining the code universe $\C_f$, take
a function that exactly performs the above described set of tasks, e.g.\
three NOTS etc.

As the subspace $A$ we simply take the generated avida code in  
\texttt{153-aagxs}.

\subsection{Comparison codes} \label{sec.comp}
To get a comparison environment, let us construct (or simulate) two,
man--made codes that perform the operations in $f$.

We use a very simple approach and simulate two different codes in $B$.
One with no loop, except for the self--copying loop, and one with
as many loops as there are logical tasks to be done. Furthermore,
we don't actually write the codes but simulate the writing using
a very rigid approach of making all the logical statements from
combinations of NANDs\footnote{This is the base for the  default reward list
for completed tasks in time slicing; see the file task.set in the Avida
package.},  
and then in detail study the needed operations just in one NAND.

This approach will give us a list of operations, and operands needed to
fulfill $f$. We can now try to look for style.

\subsection{The measure}
To start with something, we used the Halstead's measures as our
$\mufat$. That is, let 
\[\mufat=(\mbox{vocabulary, length, difficulty, volume, effort}).\]

\subsection{Simulations}
The actual simulations is done by extracting creatures, i.e.\ individual
programs, from a run of the Avida program, as described above. The
extraction consists of some data, the table of performance, but also the
program itself. We then use a $C++$ program to automatically analyze the
extraction by first reading of the data and the performance table, and then
find the Halstead's measure vector $\mufat$ for the actual code.  The $C++$
program then also uses the performance table to simulate the two different
comparison codes described in Section \ref{sec.comp} above and to calculate
their Halstead's measures respectively.  The program then also computes
$\wfat$ when $A$ has only the code of extracted individual, and $B$ consist
of the two comparison codes.  ( $\wfat$ for other combinations of $A$ and
$B$ are also considered.) The results are then exported to a Maple file
where one can more easily work with the output.


\subsection{Preliminary results}
Since $\#A$ is one, and $\#B$ is two, we can not 
talk about results in any statistical sense of course. Nevertheless,
we can present some outcome that should be seen as a 
indicator of what eventually can be done. Even if we have
a very small $\C_f$, since essentially every creature represents a
unique list of performed tasks, which leads to essentially a unique
function $f$; we have many classes $\C_f$, and we are free to 
change parameters, such as the initial random seed, the
reward table for the tasks, etc. We can therefore compare
profiles; see for example Figure~\ref{fig.8500}.

Let us more in detail see what kind of output you can get by
going back to our old friend, born after 19441 generations, \texttt{153-aagxs}.

After feeding the file \texttt{153-aagxs} into the C++ program,
which disects the code, and simulates the two comparison codes,
we get as outputs things like:
The Halstead's measures:
\[n_1 =19, 
 n_2 =3,
 N_1 =153,
 N_2 =31\]
 \[\mbox{difficulty} =98.1667\]
 \[\mbox{volume} =820.535\]
 \[\mbox{effort} =80549.2.\]
And we also get 
\[\theta=0.0085549
 \mbox{ and } \eta=391865,\]
which does not really tell you much in this meager situation ($\eta$
is huge since the variance of $X$ is so extremely tiny for this special 
case).
We can view $\wfat^+$ in a diagram; see Figure \ref{fig.w2.153}.
(Remember that we normalized the measures into $[0,1]$ by the 
transformation $\frac{x}{1+x}$ in order to make the weights in $\wfat$
play in the same division.)

\vskip 4mm

\begin{figure}[htb]
\centerline{\psfig{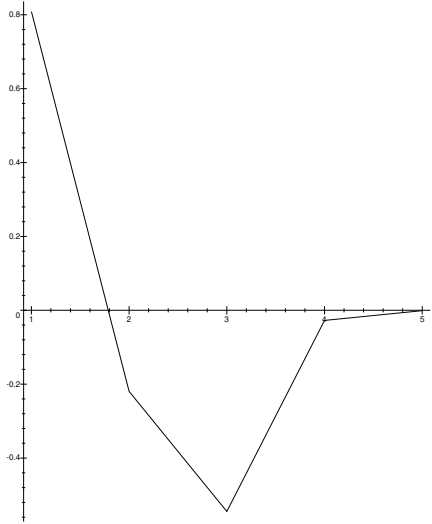}}
\textsf{\caption{\textit{Here is an example of a fingerprint after about 1000
generations using Halstead's complexity measures as $\mufat$}.  \label{fig.aik}}
}\end{figure}

\begin{figure}[htb]
\centerline{\psfig{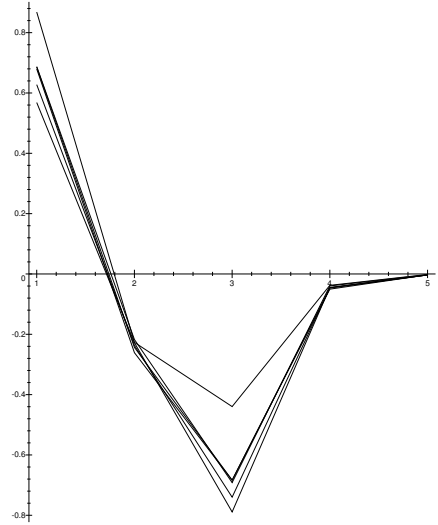}}
\textsf{\caption{\textit{Here are seven fingerprints after about 8500.
Can we hope for some convergence? And if so, what will that tell us?}  
\label{fig.8500}}
}\end{figure}

\begin{figure}[htb]
\centerline{\psfig{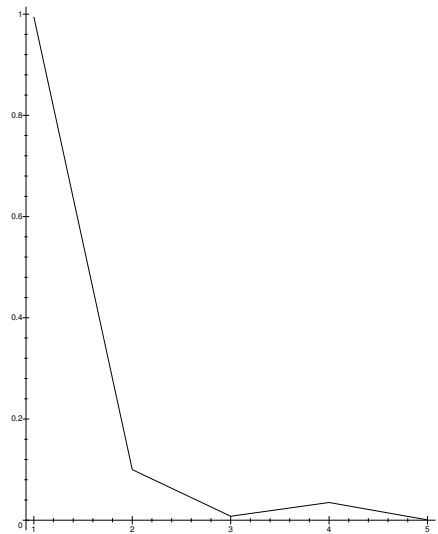}}
\textsf{\caption{\textit{Here is    $\wfat^+$ when $A$ is the single Avida
generated code} \texttt{153-aagxs} \textit{and 
$B$ consists of the no-loop-code and the all-loop-code, both
based entirely on NAND combinations.}
\label{fig.w2.153}}
}\end{figure}

\begin{figure}[htb]
\centerline{\psfig{figure=w0.153-aagxs.epsi,width=.7\hsize}}
\textsf{\caption{\textit{$\wfat^+$ when $A$ is} \texttt{153-aagxs}\textit{  and 
$B$ is the single no-loop-code.}
\label{fig.w0.153}}
}\end{figure}

\begin{figure}[htb]
\centerline{\psfig{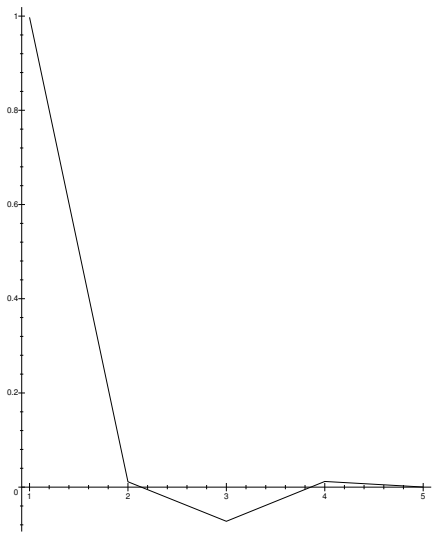}}
\textsf{\caption{\textit{Here is  $A$} \texttt{153-aagxs} \textit{and 
$B$ is the all-loop-code.}
\label{fig.w1.153}}
}\end{figure}

\begin{figure}[htb]
\centerline{\psfig{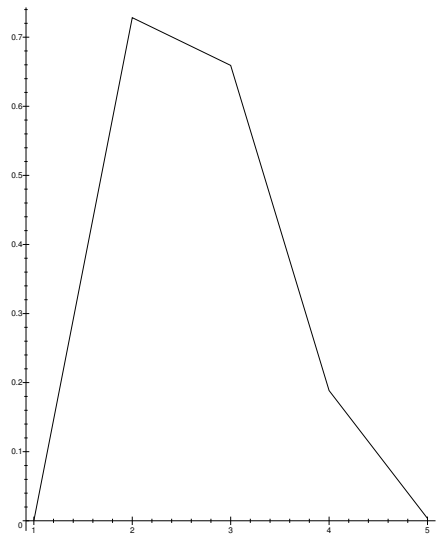}}
\textsf{\caption{\textit{$\wfat^+$ when $A$ is the all-loop-code  and 
$B$ is the  no-loop-code.} 
\label{fig.w3.153}}
}\end{figure}

\begin{figure}[htb]
\centerline{\psfig{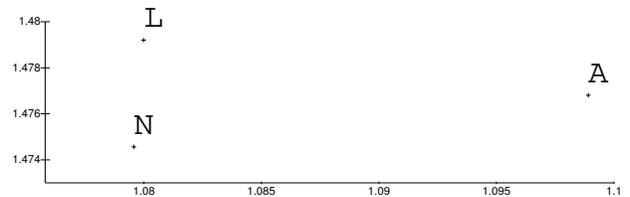}}
\textsf{\caption{\textit{Here is a picture of a principal component
analysis, see Section \ref{sec.common} of
the $\mufat$ vector for the three different codes connected to}
 \texttt{153-aagxs}. \textit{The letter}  \texttt{A} \textit{stands for the
Avida code and the letters}  \texttt{N} \textit{and} \texttt{L} 
\textit{for the no--loop code and the loop code. We see that the two
simulated codes,} \texttt{N} \textit{and} \texttt{L}, \textit{
are more together indicating that they are more similar to each other
than the evolutionary generated code from Avida.}
\label{fig.principal}}
}\end{figure}

\onecolumn



\end{document}